\documentclass[12pt]{article}
\usepackage{epsf}
\usepackage{amsmath}
\usepackage{mathptm}
\usepackage{pst-plot}
\usepackage{pst-grad}
\usepackage{pst-node}
\usepackage {amsmath}
\input amssym.def
\input amssym.tex

\def \C{{\Bbb C}}
\def \Q{{\Bbb Q}}

\def \l{{\lambda}}

\def \Z{{\Bbb Z}}

\def \O{{\cal O}}

\def \re {{\it Remark\ \ }}
\def \proof{{\noindent{\it Proof.\ \ }}}
\newtheorem{Th}{THEOREM}[section]

\newtheorem{prop}[Th]{PROPOSITION}

\newtheorem{lemma}[Th]{LEMMA}
\newtheorem{cor}[Th]{COROLLARY}

\title{Equivariant Todd Classes for Toric Varieties}
\author {Jean-Luc Brylinski\thanks{ This research was supported in part by NSF grant DMS-9504522 and DMS-9803593}\\
 Department of Mathematics\\The  Pennsylvania State University \\
University Park, PA 16802, USA\\
        {\sf email: jlb@math.psu.edu } \\
Bin Zhang$^*$\\
 Department of Mathematics\\The State University of New York,\\
Stony Brook, NY 11794-3651 USA\\
        {\sf email: bzhang@math.sunysb.edu } }
\date {}
\begin{document}
\maketitle

\begin{abstract} 
For a complete toric variety, we obtain an explicit formula for the localized equivariant Todd class in terms of the combinatorial data -- the fan. This is based on the equivariant Riemann-Roch theorem and the computation of the equivariant cohomology and equivariant homology of toric varieties. 
\end{abstract}
\newpage
\section {Introduction}

A toric variety is a normal variety with an action of an algebraic torus,  which contains the torus as an open dense subset. They form a very interesting family of complex algebraic varieties, and were introduced by Demazure \cite {dm1} and Kempf-Knudsen-Mumford-Saint-Donat \cite {kkms1}. Every toric variety can be described by a set of combinatory data - a fan. The theory of toric varieties established a now classical connection between algebraic geometry and the theory of convex polytopes. In particular, the Todd classes of complete toric varieties have important consequences for numerical results on convex polytopes, as observed by Danilov. Since the seventies, several authors have established formulas for  Todd classes of toric varieties \cite {mr1}, \cite {pj1}, \cite {ggk1} \cite {cs1}, \cite {cs2}.

Brion and Vergne \cite {bv1} proved an equivariant version of Riemann-Roch theorem for complete simplicial toric varieties, and in particular they introduced the equivariant Todd classes for simplicial toric varieties with values in equivariant cohomology. We follow this direction to work on the Todd classes of general complete toric varieties. But we use equivariant homology instead of equivariant cohomology, because in the general the Todd classes take values in homology \cite {bfm1}. 

Basically the paper is divided into two parts. In the first part, we set up the equivariant Riemann-Roch theorem following the framework of Baum-Fulton-MacPherson \cite {bfm1}. For our purpose, we only consider $G$-varieties (maybe singular) over $\C$, where $G$ is a complex linear reductive algebraic group. Roughly speaking, the problem that we need to solve is to give correct definitions. To get the proper definition, the key is Totaro's approximation of $EG$ \cite {eg2}. Based on it, we can define equivariant homology \cite {eg1} and establish the equivariant Riemann-Roch theorem (Theorem \ref {RR}). Edidin and Graham have establised an version of equivariant Riemann-Roch theorem \cite {eg3}, because our result is slightly different from theirs, as we use equivariant homology instead of the equivariant Chow group, for completeness, we include this discussion here.
 
The second part is the computation. We compute the equivariant cohomology and equivariant homology for toric varieties. Based on that, we obtain a concrete formula for the equivariant Todd class (Theorem \ref {mainth}). Here the power maps of toric varieties play a very important role.

For simplicity, in the paper, we assume all cohomolgy or homology is taken over the rational number field $\Q$.

Here we would like to thank Michel Brion, Mich\`ele Vergne, William Fulton, Sylvain Cappell, Julius Shaneson for useful correspondence, Dan Edidin and William Graham for pointing out a mistake in our previous definition of equivariant Todd classes. 

\section {Totaro's approximation of $EG$}

The purpose of this section is to describe a construction \cite {eg2}, due to Totaro, of an algebro-geometric substitute for the classifying space of topological group.

For a topological group $G$, the spaces $EG$ and $BG$ are infinite dimensional, and are hard to control in general. What Totaro produces for a reductive algebraic group is a directed system of algebraic $G$-bundles $E_n \to B_n$ with the following property:
for any principal algebraic $G$-bundle $E\to X$, there is a map $X'\to X$ with the fiber isomorphic to ${\Bbb A}^m$, such that the pullback bundle $E'\to X'$ is pulled back from one of the bundles in the directed system by a map $X'\to B_n$.

For a complex linear reductive algebraic group $G$ (In fact it works for more generality), we can construct the directed system as follows:

Every object in the directed system is a pair $(V, V')$, such that $V$ is a representation of $G$, $V'$ is a non-empty open set (in \'etale topology) of $V$ with free $G$ action, $V'\to V'/G$ is a principal $G$-bundle. $(V, V')<(W, W')$ if $V$ is a subspace of $W$, $V'\subset W'$, ${\rm codim}_V(V-V')<{\rm codim}_W(W-W')$. The morphisms are just inclusions.

For some groups, there is a convenient choice of the directed system. For example, if $G=T$ is a complex torus of rank n, then we can take $\{(V^n, (V^n-\{0\}))\ | \ {\rm dim} V=l\}_l$ as a directed system, where $V^n$
has a $T\cong (\C^*)^n$ action: $(t_1, t_2, \cdots, t_n)(v_1, v_2, \cdots, v_n)=(t_1v_1, t_2v_2, \cdots, t_nv_n)$.
 
\section {Equivariant cohomology}

Let $G$ be a topological group, $X$ be a $G$-space, Borel \cite {ba1} defined the equivariant cohomology of $X$ as
$$H_G^*(X)=H^*(X\times ^GEG)$$
where $EG\to BG$ is the universal $G$-bundle.

There is an interesting property of equivariant cohomology with respect to the change of groups. 

\begin {prop}
Let $H$ be a closed subgroup of $G$, $X$ be a $H$-space, then
$$H^*_G(G\times ^H X)\cong H^*_H(X)
$$
\end{prop}

Using the directed system in Section 2, we have another description of equivariant cohomology:

\begin {lemma} For a complex Lie group $G$ (which has a directed system as in the previous section) and $X$ a complex algebraic variety on which $G$ operates algebraically,
$$H^*_G(X)=\lim _{\underset {(V, V')} \gets} H^*(X\times ^G V')$$
where the inverse limit is taken over any directed system satisfying the properties in previous section.
\end {lemma}

\proof For any $(V, V')$ in the directed system, $V'\to V'/G$ is a principal $G$-bundle, so we have a commutative diagram:
$$\begin {array}{ccc}
V'&\to EG\\
\downarrow&\downarrow \\
V'/G&\to BG 
\end{array}
$$
There is an induced map: $X\times ^G V' \to  X\times ^GEG$ and we obtain a map: $H^*_G(X) \to H^*(X\times ^G V')$, thus a map $H^*_G(X)\to\lim _{\underset {(V, V')} \gets} H^*(X\times ^G V')$.

By the properties of the directed system, we can prove this map is an isomorphism.
$\square$

For every $G$-vector bundle $E$ over $X$, we can check $E\times ^G V'\to X\times ^G V'$ and $E\times ^G EG\to X\times ^G EG$ are all $G$-vector bundles, so we can define the Chern character $ch^G(E) \in H^*_G(X)$ and $ch^V(E) \in H^*(X\times ^G V')$ by 
$$ch^G(E)=ch (E\times ^G EG)$$
and
$$ch^V(E)=ch (E\times ^G V')$$
We see that under the map $H^*_G(X) \to H^*(X\times ^G V')$, $ch^G(E)$ goes to $ch^V(E)$.

In the same way, we can define the cohomological equivariant Todd class for a $G$-vector bundle over a smooth $X$ with smooth $G$-action. For any $(V, V')$ in the directed system, $X\times ^G V'$ is smooth. So we can define the cohomological equivariant Todd class $Td^G(X)$ to be the element in $H_G^*(X)$ which maps to $Td(X\times ^GV')$ for every $(V, V')$ in the directed system.
  
\section {Equivariant homology}

Now for a complex linear reductive algebraic group $G$, we can define the equivariant homology \cite {eg1} for a $G$-variety $X$ as 

$$H_i^G(X)=H^{BM}_{i+2l-2g}(X\times ^GV'),$$
where $\{(V, V')\}$ is a directed system for $G$ and codim$_V(V-V')$ is big enough, $H_*^{BM}$ is Borel-Moore homology, $l={\rm dim}_{\C}V$, $g={\rm dim}_{\C}G$.
It is possible that $H_i^G(X) \not =0$ for negative $i$.

\begin {prop} $H^{G}_i(X)$ is independent of $(V, V')$ if codim$_V (V-V')$ is big enough. 
\end {prop}

\proof See \cite {eg1}. $\square$
\vskip 2mm
Just as for equivariant cohomology, we have the following property.
 
\begin {prop}
Let $H$ be a closed subgroup of $G$, $X$ be a $H$-space, then
$$H_*^G(G\times ^H X)\cong H_*^H(X)
$$
\end{prop}

As usual, $H^G_*(X)$ is a module over $H_G^*(X)$. In other words there is a cap product:
$\cap:  H_G^*(X)\otimes H^G_*(X)\to H^G_*(X)$ defined in the following way:

If $\alpha \in H_G^*(X)$, 
then for any $(V,V')$, the image of $\alpha$ under the map $H_G^*(X)\to H^*(X\times ^G V')$ defines a map by cap product: $H^{BM}_*(X\times ^G V') \to H^{BM}_*(X\times ^G V') \to H^{G}_*(X)$, thus all these induce a map from $H^{G}_*(X)$ to $H^{G}_*(X)$. This gives us the module structure.

Now let us define the orientation class $[X]_G \in H^{G}_{2n}(X)$ for a variety $X$ of complex dimension $n$.

For any object $(V, V')$ in the directed system, $X\times ^G V'$ defines a class $[X\times ^G V'] \in H^{BM}_{2n+2l-2g}(X\times ^GV')$, when codim$_V(V-V')$ is big enough, this is our orientation class $[X]_G$.

\section {Equivariant Riemann-Roch Theorem}

Now we can define the equivariant Todd class $\tau ^G$ as follows \cite {eg3}, For any equivariant coherent $G$-sheaf ${\cal F}$ over $X$ ($X$ is a complex algebraic variety equipped with an algebraic $G$-action), if $p: \ X\times V' \to X$ is the projection, then the pullback sheaf $p^*{\cal F}$ is a $G$-equivariant coherent sheaf over $X\times V'$, it descents to a coherent sheaf ${\cal F}_V$ over $X\times ^G V'$. We have the Todd class 
$\tau ({\cal F}_V)$ is in $H^{BM}_*(X\times ^G V')$. Meanwhile, we have a vector bundle over $X\times ^G V'\to X\times ^G (V'\times V)$, its cohomological Todd class $Td(X\times ^G (V'\times V))$ is in $H^*(X\times ^G V')$ and is invertible. Thus
we get a map $K^G_0(X)\to H^{BM}_*(X\times ^G V')$, ${\cal F} \mapsto \tau ({\cal F}_V)\cap Td^{-1}(X\times ^G (V'\times V))$. We can check, it is compactible with the transition maps, so we get a map  $K^G_0(X)\to \Pi _p H^{G}_p(X)=\hat H^G_*(X)$, that is $\tau ^G$. 

\begin {Th} \label {RR} [Equivariant Riemann-Roch Theorem]
 $\tau ^G$ is a natural transformation (with respect to $G$-equivariant proper maps) from $K_0^G(X)$ to $\hat H^G_*(X)$, such that for any $G$-variety $X$, we have the following commutative diagram:

$$\begin {array}{cclcc} 
K_G^0(X)&\otimes& K_0^G(X) &{\overset {\otimes} \to}& K_0^G(X)\\
&\downarrow&ch^G\otimes \tau ^G&&\downarrow\tau ^G\\
\hat H_G^*(X)&\otimes& \hat H_*^G(X) &{\overset {\cap} \to} &\hat H_*^G(X)
\end {array}
$$

and if $f: X\to Y$ is a proper $G$-equivariant algebraic map, then we have following commutative diagram
$$\begin {array}{rrr} 
K_0^G(X) &{\overset {f_*} \to}& K_0^G(Y)\\
\downarrow \tau ^G&&\downarrow\tau ^G\\
\hat H_*^G(X) &{\overset {f_*} \to} &\hat H_*^G(X)
\end {array}
$$

Moreover, if $X$ is non-singular, and ${\cal O}_X$ is the structure sheaf, then 
$$
\tau ^G({\cal O}_X)=Td^G(X) \cap [X]_G
$$

Furthermore if any morphism $\phi: K_0^G(X) \to \hat H^G_*(X)$ satisfies the above properties and for any free $G$-space $X$, $\phi (X)=\tau (X/G)$, the normalized Baum-Fulton-MacPherson Todd class of $X/G$, then $\phi=\tau ^G$.
   
\end {Th}
\proof We can directly check the properties of $\tau ^G$ by definition. As the uniqueness, it follows from the fact that any space $V'$ is a free $G$-space, so is $X\times V'$, thus $\phi $ agrees with $\tau ^G$ on $X\times ^G V'$, passing to the limit, we know $\phi=\tau ^G$. $\square$

\section {Equivariant cohomology of toric varieties}

Having defined the equivariant Todd class, now let us turn to toric varieties and describe the equivariant Todd classes concretely. Before further discussion, let us recall some general facts about toric varieties.

Let $T\cong (\C ^*)^d$ be an algebraic torus of dimension $d$, $N \cong Hom (\C ^*, T)\cong \Z ^d$ be the set of $1$-parameter subgroups of $T$, $\Sigma$ be a fan, (ie, a family of strongly convex rational cones in $N \otimes \Q$ such that the intersection of any two elements is again an element in the family and each face of a cone
of the family is in the family), $X$ be the toric variety associated to it. If we denote the dual lattice of $N$ by $M$, then $M$ is isomorphic to the character group of the torus $T$. 

For any $\sigma \in \Sigma$, let us introduce some notation related to it:
$$\begin {array}{lcl}
U_{\sigma}&:& {\rm -the\ affine\ toric\ variety\ corresponging\ to\ it} \\ 
O_{\sigma}&:& {\rm -the\ corresponding\ }T-{\rm orbit\ (which \ is \ closed \ in \ U_{\sigma})}\\
\sigma ^{\perp}&:&=\{u \in M_{\Q} \ |\  u|_{\sigma}=0\}\\
\check \sigma&:&=\{u \in M_{\Q} \ |\  <u,v>\ge 0, \ v\in \sigma\}\\
\bar \sigma&:&=\check \sigma \cap M \\
S_{\sigma}&:& {\rm -the\ symmetric\ algebra\ over} of M/\sigma ^{\perp}\cap M \ {\rm over\ } \Q\\
T_{\sigma}&:& {\rm -the\ subgroup\ of}\ T\ {\rm with\ character\ group}\ M/\sigma ^{\perp}\cap M  
\end{array}
$$
Notice that the group ring $\C [\bar \sigma]$ is the regular function ring of $U_{\sigma}$ and $U_{\sigma}\to O_{\sigma}$ is a $T$-equivariant deformation retraction. In fact $U_{\sigma }\cong T\times ^{T_{\sigma}}F_{\sigma}$, where $F_{\sigma }$ is an affine $T_{\sigma}$-toric variety, more explicitly, $\sigma $ spans a linear subspace $N_{\sigma}$ of $N\otimes \Q$, in which we have a natural fan consists of $\sigma$ and its faces, $F_{\sigma}$ is the corresponding toric variety.

Let us use $\Sigma '$ to denote the subset of maximal cones in $\Sigma$. $\Sigma$ has a partition $\{\Sigma (i)\}$ according to the dimension of cones, $\Sigma (i)=\{\sigma \in \Sigma\ |\ {\rm dim}\sigma =i\}$.

For a $T$-space $Y$, following Borel \cite {ba1}, $H_T^*(Y)$ by definition is $H^*(Y\times ^TET)$, where $ET\to BT$ is the universal $T$-bundle. When $Y$ is a point, we have $H^*_T({\rm point})=H^*(BT)\cong \Q[x_1, x_2, \cdots, x_d]$.

\begin {lemma}
$$H^*_T(U_{\sigma})=S_{\sigma}$$
\end {lemma}

\proof We already know that, as a $T$-space, $U_{\sigma}$ is $T$-homotopic to $O_{\sigma}\cong T/T_{\sigma}$, so $H^*_T(U_{\sigma})\cong H^*_T(T/T_{\sigma})\cong H^*_{T_{\sigma}}({\rm point})\cong S_{\sigma}$.$\square$

\begin {lemma} If $\tau \subset \sigma$ is a face, we have the inclusion
$i: U_{\tau} \to U_{\sigma}$ \cite {fw1}, then $i^*: H^*_T(U_{\sigma}) \to H^*_T(U_{\tau})$ is the map $S_{\sigma} \to S_{\tau}$ induced by $M\otimes \Q/\sigma ^{\perp} \to M\otimes \Q/\tau ^{\perp}$. 
\end {lemma}

\proof See \cite {bv1}. $\square$

For the toric variety $X$ associated with the given fan $\Sigma$, notice that $X\times ^T ET$ has an open covering ${\cal C}$: $\{U_{\sigma}\times ^T ET\}_{\sigma\in \Sigma '}$, thus we obtain a spectral sequence which converges to $H^*_T(X)$ \cite {hf1} with

$$E_1^{p,q}=\oplus _{\{\sigma_0, \sigma_1, \cdots, \sigma_p\}\subset \Sigma '}H^{q}(U_{\sigma_0\cap \sigma_1 \cap \cdots \cap\sigma _p}\times ^T ET )$$

\begin{Th}
This spectral sequence degenerates at $E_2$, and  we have a canonical isomorphism of graded $H^*(BT)$ algebras,
$$Tot(E_2)\cong H_T^*(X)$$
\end{Th}

We will prove this theorem in next section.
Therefore, every element in $H^n_T(X)$ can be represented as a class $(a_{\sigma_0, \sigma_1, \cdots, \sigma_p,q})$ with  $a_{\sigma_0, \sigma_1, \cdots, \sigma_p,q}\in H^{q}(U_{\sigma_0\cap \sigma_1 \cap \cdots \cap\sigma _p}\times ^T ET )$ and $p+q=n$.
\vskip 3mm
In case $X$ is simplicial, notice that for any $q$, the sequence
$$0 \to E_1^{0,q}\to E_1^{1,q} \to E_1^{2, q}\to \cdots $$
is exact except at $E_1^{0,q}$, we see 

\begin {Th} \cite {bv1}
If $X$ is simplicial, then $H^*_T(X)$ is isomorphic to the algebra of continuous piecewise polynomial functions. 
\end{Th}

\section {The power maps of toric varieties}

This is a very special property of toric varieties. For any positive integer $n$, we have a map:
$$\begin {array}{cccc}
\l _n: &\Sigma& \to &\Sigma\\
&x & \mapsto & nx
\end{array}
$$
It induces a map \cite {fw1}: $X \to X$, which we again denote by $\l _n$. If we restrict $\l _n $ to the open dense orbit $T$, then it is the map $T\to T: \ t\mapsto t^n$, so we call $\l _n$ the $n-th$ power map.
$\l _n$ is not $T$-equivariant in general (except for $n=1$). Let us first consider some properties of these power maps.

\begin{lemma}
For a continuous group homomorphism $\phi: G_1 \to G_2$ of topological groups, $\phi$ induces  a $G_1$-equivariant map
$EG_1 \to EG_2$ and a commutative diagram:
$$\begin {array}{ccc}
EG_1 &{\overset \phi \to}&EG_2 \\
\downarrow &&\downarrow \\
BG_1 &{\overset \phi \to}& BG_2 
\end{array}
$$
\end{lemma}

X has a canonical $T$-action. Based on this action, we can construct a family of $T$-actions $\{\psi _n\}$on $X$:
$$\begin {array}{cccc}
\psi _n: &T\times X& \to &X\\
&(t,x) & \mapsto & t^n.x
\end{array}
$$
We will denote by $_nX$ the space $X$, on which $T$ acts by $\psi _n$. Then $\l _n: X \to _nX$ is an equivariant $T$-map.

Apply the above lemma to the $n$-th power endomorphism of $T$:
$\l _n: T\to T$ we get a commutative diagram of principal $T$-fibration:
$$\begin {array}{ccc}
ET&{\overset {\l _n} \to}&_nET\\
\downarrow&&\downarrow \\
BT&{\overset {\l _n} \to}&BT
\end{array}
$$

\begin{lemma}
$\l _n: BT\to BT$ induces an isomorphism $\l  _n ^*: H^*(BT) \to H^*(BT)$, and $\l _n ^*$ acts on $H^{2q}(BT)$ by $n^q$.
\end{lemma}

\proof When $d=1$, we have a concrete model for $BT$: $\C {\Bbb P} ^{\infty}$ and the map $\l _n$ is nothing but $[x_0: x_1: x_2\cdots ]\mapsto [x_0^n: x_1^n: x_2^n\cdots ]$. For general $d$, we can take $BT\cong (\C {\Bbb P} ^{\infty})^d$, and in a similar way, we can use the K\"unneth theorem to prove this lemma. $\square$

\begin {lemma} For any $T$-space $X$, we have 
$$
H^*_T(X)\cong H^*(X\times ^T \ _nET)$$
\end{lemma}
\proof We have a Cartesian diagram of fibrations with fiber $X$:
$$\begin {array}{ccc}
X\times ^T ET&{\overset {id \otimes \l _n} \to}&X\times ^T\ _nET\\
\downarrow&&\downarrow \\
BT&{\overset {\l _n} \to}&BT
\end{array}
$$

So we get a map between Leray spectral sequences; because $BT$ is simply connected, the $E_2$ terms are given by tensor product. This implies that the map between Leray spectral sequences yields an isomorphism at the $E_2$ level. This proves our claim. $\square$

\begin{lemma}
Under the isomorphism $H^*_T(U_{\sigma}) \cong S_{\sigma}\cong \Q[x_1, x_2, \cdots, x_{{\rm dim}\sigma}]$, the induced map $(id \otimes \l _n)^*: H^*_T(U_{\sigma}) \to H^*_T(U_{\sigma}\times ^T \ _nET)$ sends $x_i$ to $n x_i$.
\end{lemma}

Now we can prove Theorem 6.3.
\vskip 3mm

\noindent {\it Proof of Theorem 6.3} Consider the map: $id \otimes \l _n: X\times ^T ET \to X \times ^T  \ _nET$, the open covering ${\cal C}$ of $X\times ^TET$ and the covering ${\cal C}'=\{U_{\sigma}\times ^T \ _nET\}$ of  $X \times ^T  \ _nET$. The map $id \otimes \l _n$ is compatible with the coverings ${\cal C}$ and ${\cal C}'$, so we get a map between the spectral sequences for $X\times ^T ET$ and $X\times ^T\ _nET$, which we constructed in previous section.

So we have the diagram:
$$\begin {array}{ccc}
E_r^{p,q}&{\overset {d_r}\to}&E_r^{p+r, q-r+1}\\
\downarrow &&\downarrow \\
E_r^{p,q}&{\overset {d_r}\to}&E_r^{p+r, q-r+1}
\end{array}
$$
where the downarrow maps are the map $(id \otimes \l _n)^*$. If we identify $E_r^{p,q}$ as a quotient of a subset of $ E_1^{p,q}$, the map $(id \otimes \l _n)^*$ is multiplication by $n^ {[q/2]}$. So we see, $d_r =0$ for $r>1$.

Because this spectral sequence degenerates at $E_2$, for any $p,\ q$, we have a filtration of $H^{p+q}_T(X)$:
$$
H^{p+q}_T(X)= F^0H^{p+q}_T(X)\supset F^1H^{p+q}_T(X)\cdots \supset F^kH^{p+q}_T(X) \cdots
$$
and the following exact sequence:
$$
0\to F^{p+1}H^{p+q}_T(X) \to F^pH^{p+q}_T(X)\to E^{p,q}_2\to 0
$$
Now using the maps $(id \otimes \l _n)_*$, we obtain the isomorphism of $H^*(BT)$ algebras:
$$Tot(E_2)\cong H_T^*(X)
$$
$\square$

\section {Equivariant homology for toric varieties}

The equivariant homology is defined as $H^T_i(X)=\ H^{BM}_{i+2(k-1)d}(X\times ^TV_k)$, where $V_k=(\C ^k-\{0\})^d$ with the diagonal action of $T$, $(t_1, t_2, \cdots, t_d).(v_1, v_2, \cdots, v_d)$ = $(t_1 v_1, t_2 v_2, \cdots, t_d v_d)$. As in section 7, $_nV_k$ means $V_k$ with the action of $n$-th power of $T$.

\begin {lemma} For any $\sigma \in \Sigma $, 
$$H_*^T(O_{\sigma})\cong S_{\sigma}[O_{\sigma}]$$
In this case, for any $a \in M/\sigma ^{\perp}\cap M$, deg $a=-2$, and $[O_{\sigma}]$ is the orientation class, deg $[O_{\sigma}]$ =2codim $\sigma$.
\end {lemma}
\proof We know $O_{\sigma}\cong T/T_{\sigma}$, so the statement reduces to the case of one point by Proposition 4.2. $\square$

For any $T$-space $X$, we have a map $id \otimes \l _n: X\times ^TV_k \to X\times ^T \ _n V _k$, 
so similarly to the proof of Lemma 7.3, we can prove

\begin{lemma} For any $T$-space $X$, 
$${\underset {\overset \to k}{\lim}}H^{BM}_*(X\times ^TV_k)\cong 
{\underset {\overset \to k}{\lim}}H^{BM}_*(X\times ^T\ _nV_k)
$$
In other words, we can use either $\{V_k\}$ or $\{_nV_k\}$ to define the equivariant homology.
\end{lemma}

Let us denote the map ${\underset {\overset \to k}{\lim}}H^{BM}_*(X\times ^TV_k)\cong 
{\underset {\overset \to k}{\lim}}H^{BM}_*(X\times ^T\ _nV_k)
$ by $(id\otimes \l _n)_*$

\begin {lemma}
Under the isomorphism $H^T_*(O_{\sigma}) \cong S_{\sigma}[2 {\rm codim}\sigma]$, the map:  $(id \otimes \l _n)_*$ maps $a \in M/\sigma ^{\perp}\cap M $ to $n a$, i.e.,  it acts on $H^T_{2k}(O_{\sigma})$ as multiplication by $n^{{\rm codim} \sigma -k}$ .
\end{lemma}
\proof It is similar to the proof of Lemma 7.2. $\square$

Let $X_k=\cup _{{\rm dim}\sigma \ge d-k} O_{\sigma}$, then $X_k$ is closed and  we obtain a $T$-invariant filtration
of $X$,
$$\{X_0 \subset X_1\subset \cdots \subset X_n=X\}$$

\begin{lemma}
If $Z\subset X$ is a closed $T$-subset of $X$, then we have following long exact sequence:
$$\cdots \to H^T_p(Z) \to  H^T_p(X)\to H^T_p(X-Z)\to H^T_{p-1}(Z)\to \cdots$$
\end{lemma}
\proof This comes from the facts that we have an exact sequence for the pair $(X\times ^T V_k, Z\times ^T V_k)$ and taking direct limits is an exact funtor. $\square$
 
\begin {lemma} For any $p$ and $k$, we have $H^T_{2k+1}(X_p)=0$
\end{lemma}

\proof Notice $X_0$ is a finite set, and for any $0 \le p \le d$, $X_p-X_{p-1}$ is the disjoint union of $p$-dimensional orbits, so the result follows by induction. $\square$

\begin{Th} \label {th-equ-h} For any $k$, we have a canonical decomposition 
$$H^T_{2k}(X)\cong \oplus _{\sigma} H^T_{2k}(O_{\sigma})$$
Therefore, every element in $H^T_n(X)$ can be represented as $(b_{\sigma})$ with $b_{\sigma} \in H^T_n(O_{\sigma})$.

\end{Th}

\proof By the above lemmas, for any $p$ and $k$, we have an exact sequence:
$$0\to H^T_{2k}(X_{p-1})\to H^T_{2k}(X_{p})\to H^T_{2k}(X_p-X_{p-1})\to 0$$
Now considering the map $(id \otimes \l _n)_*$, we have the following commutative diagram:
$$\begin {array}{ccccccccc}
0&\to &H^T_{2k}(X_{p-1})&\to& H^T_{2k}(X_{p})&\to &H^T_{2k}(X_p- X_{p-1})&\to& 0\\
&&\downarrow&&\downarrow&&\downarrow \\
0&\to &H^T_{2k}(X_{p-1})&\to &H^T_{2k}(X_{p})&\to &H^T_{2k}(X_p- X_{p-1})&\to & 0
\end{array}
$$
where the downarrow maps are the maps $(id\otimes \l _n)_*$. If $i: \ X_p-X_{p-1} \to X_p$ is the inclusion, for any $a \in H_{2k}^T(X_p)$, then $i_*(id\otimes \l _n)_*(a)=n^{p-k}i_*(a)$, so 
$(id\otimes \l _n)_*(a)=n^{p-k}a+b$, where $b\in H^T_{2k}(X_{p-1})$. It is easy to see there exists $c \in H^T_{2k}(X_{p-1})$ such that $(id\otimes \l _n)_*(c)=n^{p-k}c+b$, that implies $(id\otimes \l _n)_*(a+c)=n^{p-k}(a+c)$. Now it is easy to see that $Ker ((id\otimes \l _n)_*-n^{p-k})\cong H^T_{2k}(X_p, X_{p-1})$ and the short exact sequence splits. So by induction, we know $H^T_{2k}(X)\cong \oplus _{\sigma} H^T_{2k}(O_{\sigma})$ and the map $(id\otimes \l _n)_*$ is scalar multiplication by $n^{{\rm codim} \sigma -k}$ on the summand $H^T_{2k}(O _{\sigma}$.
$\square$

So far, we explored the structure of $H^T_*(X)$ as an Abelian group. As usual, $H^T_*(X)$ is a $H^*_T(X)$ thus a $H^*(BT)$-module, so let us consider the module structure now.

\begin {lemma}
\label {cor-f-class} The homology class $[\bar O_{\sigma}]\in H^T_*(X)$ under the canonical isomorphism is $(a_{\tau})$, where $a_{\tau}=1$ if $\tau=\sigma$, $=0$ otherwise.
\end{lemma}
\proof Let $p=codim \ \sigma$, then $[\bar O_{\sigma}]\in H^T_{2p} (X_p)\subset H^T_{2p}(X)$ is of degree $2p$, and the exact sequence we used in Theorem \ref {th-equ-h} is
$$0\to H^T_{2p}(X_{p})\to H^T_{2p}(X_p- X_{p-1})\to 0$$
so the lemma is clear. $\square$
\vskip 2mm

$H^T_*(X)$ is generated as an $H^*(BT)$-module by the classes $\{[\bar O_\sigma ]\}$, we have relations of the following type.

\begin {lemma} \label {relation} For any $\sigma \in \Sigma$, $\l \in \sigma ^{\perp}$ we have 
$$\l \cap [\bar O_{\sigma}]=\sum _{\tau} <\l , n_{\sigma  \tau}>[\bar O_{\tau}]$$
where the summation is taken over all these $\tau $'s containing $\sigma $ as a facet, $ n_{\sigma  \tau} \in N/(\sigma \cap N)$ is the unique generator of the semigroup $(\tau \cap N)/ (\sigma \cap N)$.
\end{lemma}

So, combining this lemma with Theorem \ref {th-equ-h}, we see that $H^T_*(X)$ is the $H^*(BT)$-module, generated by $\{[\bar O_\sigma]\ | \ \sigma \in \Sigma\}$, modulo the relations given in Lemma \ref {relation}. It is isomorphic to the equivariant Chow group \cite {bm1}.  

An immediate consequence of this Lemma \ref {relation} is the localization theorem for equivariant homology.

\begin {cor} (Localization theorem) For a toric variety $X$, $i: X^T \to X$ induces an isomorphism of equivariant homologies up to $H^*(BT)$ torsion.
\end {cor}

In the case of complete simplicial fans, 
this description of equivariant homology is dual to the equivariant cohomology. Assuming $X$ is complete and simplicial, Brion-Vergne proved in \cite {bv1} that $H^*_T(X)$ is isomorphic to the algebra $R_{\Sigma}$ of continuous piecewise polynomial functions on $\Sigma$. We will write down the Poincare duality isomorphism  $PD$ between $R_{\Sigma}$ and $H_*^T(X)$ explicitly.

\begin {lemma} \label {isom-p-h} For a complete simplicial toric variety $X$ of complex dimension d associated to the fan $\Sigma$, we have an $H^*(BT)$-module morphism $PD: R^*_{\Sigma}\cong H_{*}^T(X)$, such that for any $k \in \Z$ 
 $$PD: R^k_{\Sigma}\cong H_{2d-2k}^T(X)$$
where $ R^k_{\Sigma}$ is the space of continuous piecewise polynomial functions of homogeneous degree $k$.
\end{lemma}
  
\proof For $\sigma \in \Sigma$, the subgroup of $N$ generated by the intersections of $N$ and all edges of $\sigma$ is a subgroup of finite index in the subgroup of $N$ generated by $N\cap \sigma$, let us denote the index by $mult (\sigma)$. 

For any $\tau \in \Sigma (1)$, there is a unique piecewise polynomial function $\xi _{\tau}$ on $\Sigma$, the so called Courant function associated to the edge $\tau$ \cite {bl1}. In fact, $\xi _{\tau}$ is a piecewise linear function, which can be described as follows. If $\sigma$ is a maximal cone, spanned by $\tau, v_2, \cdots, v_d$, and if $\{\tau ^*, v_2 ^*, \cdots, v_d ^*\}$ is the dual basis of $M\otimes \Q$, then $\xi _{\tau} |_{\sigma}=\tau ^*$. It is clear that $\xi _{\tau} |_{\sigma}=0$ if $\sigma$ does not contain $\tau$. 

Thus for $\sigma \in \Sigma$, we have a continuous piecewise polynomial function $\phi _{\sigma}$ on $\Sigma$, of degree 2dim$\sigma$:
$$\phi _{\sigma}=mult(\sigma) \Pi _{\tau \in \sigma (1)} \xi _\tau$$
$\phi _{\sigma}$ vanishes identically on all cones which do not contain $\sigma$.

Now we can define two $H^*(BT)$-module homomorphisms. By \cite {bl1}, $R_\Sigma$ is a $S(M)$-module with generators $\{\phi _\sigma \}$, where $S(M)$ is the symmetric algebra of $M$ over $\Q$. So we can define a map $F: R_\Sigma \to H^T_*(X): $ $\sum a_\sigma (v)\phi _\sigma \mapsto \sum (-1)^{dim \sigma} a_\sigma (-v)[\bar O _\sigma]$. Obviously it is a $H^*(BT)$-module homomorphism. Also we can define a natural homomorphism: $G: H^T_*(X)\to R_\Sigma$, which is given by $\sum a_\sigma (v)[\bar O_\sigma] \mapsto \sum (-1)^{dim \sigma} a_\sigma (-v)\phi _\sigma$.

For a simplicial toric variety, as $H^*(BT)$-modules, $H^T_*(X)$ and $H^T_*(X)$ have no torsion. If we take the multiplicative set $L\subset H^*(BT)$ to be $\{\Pi m_i \ | \ m_i \in M-\{0\}\}$, then both $H^*_T(X) \to L^{-1}H^*_T(X)$ and $H_*^T(X) \to L^{-1} H_*^T(X)$ are embeddings. 

We can describe the localized maps $L^{-1}F$ and $L^{-1}G$ as follows.
Any $f \in R_{\Sigma}$ is determined by $\{f|_{\sigma}\}_{\sigma \in \Sigma '}$, where $f|_{\sigma}$ is the restriction of $f$ to the maximal cone $\sigma$. The localized  map $L^{-1}F: L^{-1} R_{\Sigma}\to L^{-1} H^T_*(X)$ is given by 
$$L^{-1}F(f)=\sum_{\sigma \in \Sigma '} \frac {f|_\sigma}{\phi _\sigma}[\bar O_\sigma]
$$
It is an isomorphism with inverse map $L^{-1}G$.

So from the following commutative diagram:

$$\begin {array}{ccc} 
R_{\sigma} &{\overset {F} \to}& H^T_*(X)\\
\downarrow &&\downarrow\\
L^{-1}R_\Sigma&{\overset {L^{-1}F} \to} &L^{-1}H^T_*(X)
\end {array}
$$
we see $F$ is a module isomorphism (with inverse $G$). This is the Poincare duality isomorphism, let us use $PD$ to denote it. $\square$

\re From the proof we see, for an element in $L^{-1} H^T_*(X)$, if its image under $L^{-1}PD^{-1}$ is in $R_\Sigma$, then it is in fact in the image of $H_*^T(X) \to L^{-1} H^T_*(X)$, thus it is an element of $H_*^T(X)$. We will use this result in the next section.
 
For a general fan, the associated toric variety may be singular. In the category of toric varieties, we have an easy way to resolve the singularities - subdividing the fan. The equivariant homologies behave as follows under subdivision:

\begin {lemma} If $\ \Sigma _1$ is a subdivision of $\ \Sigma$, then the induced map of $f: X_{\Sigma _1}\to X_{\Sigma}$  is the following: $f_*(\sum a_{\sigma '}[\bar O_{\sigma '}])=\sum b_{\sigma}[\bar O_\sigma ]$, with $b_{\sigma }=\sum a_{\sigma '}$, where we sum over all $\sigma '\in \Sigma _1$ such that: $ \sigma '\subset \sigma$ and dim $\sigma '$=dim $\sigma$. 
\end{lemma} 
\proof Notice that $f$ induces an isomorphism $O_{\sigma '}\cong O_{\sigma}$ and $f_*$ is an $H^*(BT)$-homomorphism. $\square$

\section {Equivariant Todd classes}

Let $M$ be a lattice of rank $d$, we will use $\Q[M]$ to denote the group ring of $M$ over $\Q$ and $S(M)$ to denote the symmetric algebra of $M$ over $\Q$. If $M=\oplus \Z e_i$, then $\Q [M]=\{\sum a_m e^m \ | \ m\in M\}$ and $S(M)\cong\Q [x_1, x_2, \cdots, x_d]$, the polynomial ring of $d$ indetermints. Let $\hat S(M)$ be the completion of $S(M)$, ie, $\Q[[x_1, x_2, \cdots, x_d]]$, We have a  map $\Q [M]\to \hat S(M)$ which is defined by $e^{x_i}\to \sum \frac {x_i^n}{n!}$ on generators and extending it to an algebra homomorphism.

Let $\Q[[M]]$ be the set of all formal power series $\sum _{m\in M} a_m e^{m}$ with rational coefficients. 
Following Brion-Vergne \cite {bv1}, we call $f\in \Q[[M]]$ summable if there exist $P\in \Q[M]$ and a finite sequence $(m_i)_{i\in I}$, such that in $\Q[[M]]$ we have,  
$$f\Pi (1-e^{m_i})=P,$$
in the case we define the sum of $f$ as:
$$S(f)=P\Pi (1-e^{m_i})^{-1}$$
Notice $S(f)$ is an element of the fraction field of $\Q [M]$.
\vskip 2mm
By this summation, we see easily the following easy but critical result:
\begin {lemma} For any $m \in M$,$\sum _{k=-\infty}^{\infty}e^{km}$ is summable and 
$$S(\sum _{k=-\infty}^{\infty}e^{km})=0$$
\end{lemma}
\vskip 2mm
\begin {lemma} For any $\sigma \in \Sigma$, $m\in M$, denote by $U_{\sigma}\times \C m \to U_{\sigma}$ the trivial line bundle on $U_{\sigma}$ on which $T$ acts by character $m$, then we have
$$ch^T(U_{\sigma}\times \C m)=e^m$$
\end{lemma}
\vskip 2mm

The multiplicity of any character $m$ of $T$ in ${\cal O}_X(U_{\sigma})\cong \C [\check \sigma]$ is either zero or one, let us denote it by $mult_{\sigma}(m)$. Then for 

$$A_{\sigma}=\sum_{m\in \bar \sigma} mult(m)e^m.$$
we know, 

\begin {lemma} $A_{\sigma}$ is summable.
\end{lemma}
\proof See \cite {bv1}. $\square$

So we can think of $S(A_{\sigma})$ as an element of the quotient field of $\Q [M]$.
\vskip 2mm

\begin {Th} \label {mainth} For a complete toric variety $X$, we have  
$$ L^{-1}\tau ^T(X)=\sum_{\sigma \in \Sigma'} S(A_\sigma)[\bar O_{\sigma}]$$
in $L^{-1}\hat H^T_*(X)$.
\end{Th}

We will prove this theorem using the following lemmas. Before the proof, we explain what the theorem means. By definition, $S(A_{\sigma})$ is a rational function in variables $e^{m_i}$ with possible denominators of the form $\Pi (1-e^{m_i})$, so it can be viewed as a Laurent series in $m_i$, thus as an element of $L^{-1}\hat H^T_*(X)$. 

\begin {lemma} Theorem \ref {mainth} holds for any smooth toric variety $X$.
\end{lemma}
\proof Because $X$ is smooth, $\tau ^T(\O _X)=Td^T(TX)\cap [X]_T$. 

If $\sigma $ is spanned by $v_1, v_2, \cdots, v_p$, and $\{v_1, v_2, \cdots, v_d\}$ spans a maximal cone containing $\sigma$, $\{v_1^*, v_2^*, \cdots, v_d^*\}$ is the dual basis, $i: U_{\sigma} \to X$ is the inclusion map, $i^*(TX)\cong U_{\sigma} \times (\oplus_{i=1}^d \C v^*_i)$. So $Td^T(i^*(TX))$ = $\Pi^d_{i=1} \frac {v_i^*}{1-e^{-v_i^*}}$, whose image in $S_{\sigma}$ is $\Pi^p_{i=1} \frac {v_i^*}{1-e^{-v_i^*}}$. We know $i^*(Td^T(TX))=Td ^T(i^*TX)$, so as a piecewise continuous polynomial function, $Td^T(TX)$ is $\Pi^p_{i=1} \frac {v_i^*}{1-e^{-v_i^*}}$ on $\sigma$. 

On the other hand, in this case $mult_{\sigma} (m, \O _X)=1$ for $m \in \bar \sigma$. For a maximal cone $\sigma$ spanned by $\{v_1, v_2, \cdots, v_d\}$, $\bar \sigma=\Z_{\ge 0}v_1^*+\cdots +\Z_{\ge 0}v_d^*$, $S(A_{\sigma})=\Pi^d_{i=1} \frac {1}{1-e^{v_i^*}}$. 
So $\sum S(A_\sigma)[\bar O_\sigma]$ is an element of $L^{-1} \hat H^T_*(X)$.

Notice $L^{-1}PD^{-1}(\sum S(A_\sigma)[\bar O_\sigma])$ is equal to $Td^T(X) \in L^{-1}R_\Sigma$. So by the comment after Lemma 8.10, this is the equivariant Todd class. $\square$

\re We proved more in this Lemma. We proved that for a smooth toric variety (this is also true for simplicial case), $L^{-1}PD^{-1}(\sum S(A_\sigma)[\bar O_\sigma])$ belongs to $\hat H^*_T(X)$ and $PD(L^{-1}PD^{-1}(\sum S(A_\sigma)[\bar O_\sigma]))$ is the equivariant Todd class.

For a cone $\sigma$, let us denote $S(\sum _{m\in \bar \sigma}e^m)$ by $S(\C [\bar \sigma])$, then we have the following lemma.

\begin {lemma}
If $\{\sigma _1, \sigma _2, \cdots , \sigma _p\}$ is a subdivision of $\sigma$, with dim $\sigma _i$=dim $\sigma$, then we have $\sum S(\C [\bar \sigma _i])=S(\C [\bar \sigma])$, where both sides are viewed as elements of the fraction field $\Q[M/\sigma ^{\perp} \cap M]$.
\end{lemma}
\proof Fist let us prove this lemma in case $p=2$. If $\tau $ is the common facet of $\sigma _1$ and $\sigma _2$, then by \cite {fw1}, there exists $u\in \check \sigma _1 \cap \check {(-\sigma _2)}\cap M$, such that $\tau =\sigma _1 \cap u^{\perp}=(-\sigma _2)\cap u^{\perp}$, and
$\bar \tau =\bar \sigma _2+\Z_{\ge 0}u$.

In the vector space spanned by $\sigma$, for any facet $\tau $ of $\sigma$, there is a unique (up to scalar multiplication) $u_{\tau} \in \check \sigma \cap M$, such that $\tau =\sigma  \cap u^{\perp}$ and $\check \sigma $ is spanned by all these $u_{\tau}$'s. So we can assume that $\check \sigma$ is spanned by $u_1, u_2, \cdots, u_p$, $\check \sigma _1$ is spanned by $u_1, u_2, \cdots, u_s, u$ and $\check \sigma _2$ is spanned by $u_{s+1}, \cdots, u_p, -u$ with $u_i \in M$. Then we see that for any $v\in \check \sigma _1 \cup \check \sigma _2$, $v +ku \in \check \sigma _1 \cup \check \sigma _2$ for any $k \in \Z$. So $S([\C [\bar \sigma _1]]+S([\C [\bar \sigma _i]]))=S(\sum _{m \in \bar \sigma _1 \cup \bar \sigma _2} e^m)$. Notice that for any character $m$ of $T$, it appears in $\bar \sigma _i$ at most once, and $\bar \sigma _1 \cap \bar \sigma_2$=$\bar \sigma$, 
by the fact $S(\sum e^{km})=0$ for any character $m$ of $T$, we see $\sum S([\C [\bar \sigma _i]])=S([\C [\bar \sigma]])$ in this case.

For the general case, notice that for a subdivision $\{\sigma _1, \sigma _2, \cdots , \sigma _p\}$, we can subdivide them  further into $\{\tau _1, \tau _2, \cdots , \tau _q\}$, such that $\sigma _i$ is one of the $\tau _j$'s or the union of two $\tau _j$'s, $\tau _j$ and $\tau _{j+1}$ have only a common facet, and $\{\tau _1, \tau _2, \cdots , \tau _q\}$ can be divided into two sets, $\{\tau _1, \tau _2, \cdots , \tau _k\}$ and $\{\tau _{k+1},\cdots , \tau _q\}$, with the following properties:

$\tau _1\cup \tau _2\cup \cdots \cup \tau _j$ is a convex cone for any $j\le k$,

$\tau _{k+1}\cup\cdots \cup\tau _l$ is a convex cone for any $l\le q$.

So by the result about the case of subdivision into two cones, the lemma is proved. $\square$

\vskip 2mm
We now return to the proof of theorem \ref {mainth} for a general toric variety $X$. Take a suitable subdivision $\Sigma _1$ of $\Sigma $ such that the corresponding toric variety $X'$ of $\Sigma _1$ is smooth \cite {fw1}. Then the induced map    
 $f:X'\to X$ is a proper map, and $f_*\O _{X'} =\O _X$, $R^if_*\O _{X'}=0$ for $i>0$, so  $\tau ^T (\O _X)=f_*\tau ^T(\O _{X'})$. By Lemma 8.8, we know $L^{-1}\tau ^T (\O _X)$ is given by the formula. This concludes the induction step. $\square$

\end{document}